# Remarks on consistency of posterior distributions

Taeryon Choi[1] and R. V. Ramamoorthi[2]

*Inha University and Michigan State University*

**Abstract:** In recent years, the literature in the area of Bayesian asymptotics has been rapidly growing. It is increasingly important to understand the concept of posterior consistency and validate specific Bayesian methods, in terms of consistency of posterior distributions. In this paper, we build up some conceptual issues in consistency of posterior distributions, and discuss panoramic views of them by comparing various approaches to posterior consistency that have been investigated in the literature. In addition, we provide interesting results on posterior consistency that deal with non-exponential consistency, improper priors and non i.i.d. (independent but not identically distributed) observations. We describe a few examples for illustrative purposes.

## Contents

## 1. Introduction

Let $\theta$ be an unknown parameter and $X_1, X_2, \ldots, X_n$ be $n$ random variables whose joint distribution is $P_\theta^{(n)}$. In order to draw inferences on $\theta$, a Bayesian posits a prior distribution $\Pi$ for $\theta$ and updates this prior to the posterior distribution given $X_1, X_2, \ldots, X_n$, which we denote by $\Pi(\cdot | X_1, X_2, \ldots, X_n)$. This paper focuses on some issues related to an asymptotic aspect of this posterior distribution, namely, consistency.

[1]Department of Statistics, Inha University, Incheon, Korea, e-mail: trchoi@gmail.com

[2]Department of Statistics and Probability, Michigan State University, A413 Wells Hall, East Lansing, MI 48824-1027, USA, e-mail: ramamoor@stt.msu.edu

The sequence of posterior distributions $\{\Pi(\cdot|X_1, X_2, \ldots, X_n)\}$ is said to be consistent at $\theta_0$, if the posterior converges, in a suitable sense, to the degenerate measure at $\theta_0$.

Posterior consistency is a kind of frequentist validation of the updating method. If an oracle were to know the true value of the parameter, posterior consistency ensures that with enough observations one would get close to this true value. Posterior consistency also assures that as more and more observations accumulate, the observations have to dominate the role of the prior in inference. There are other interpretations related to merging of opinions and other concepts. We refer the reader to Diaconis and Freedman [10].

In order to set the perspective for this paper we begin with a short summary of earlier results in posterior consistency. Details and additional references can be found in Ghosh and Ramamoorthi [19]. The first posterior consistency result goes back to Laplace. In more recent times posterior consistency and asymptotic normality of the posterior were established for regular finite dimensional models. In a seminal paper, Freedman [12] gave a nonparameteric example, with integer-valued observations, where the posterior is inconsistent. In [10], Diaconis and Freedman showed that in the nonparametric case inconsistency can occur, even in location models with an Euclidean parameter. They suggested that instead of searching for priors that would be consistent at all unknown values of the parameter it would be fruitful to study natural priors and identify points of consistency.

On the positive side, Freedman [12, 13] and soon after Schwartz [25] provided conditions under which the posterior probability of a set $A$ will go to 0. These conditions involved two parts, one on prior positivity of Kullback–Leibler neighborhoods and the other on existence of certain test functions. Under the assumption of prior positivity of Kullback–Leibler neighborhoods, Barron gave necessary and sufficient conditions for the posterior probability of $A$ to go to 0. These results were then specialized to weak and $L_1$ neighborhoods by Barron et al. [3], Ghosal et al. [15] and Walker [32].

One aspect of these results was that they all established exponential consistency. In this paper we first give a quick review of these results from a slightly different perspective with a focus on the role of exponential consistency. We then give an example where there is consistency but not exponential consistency. The example also shows that the exponential aspect is not driven by the Kullback–Leibler condition.

Another early result in consistency is due to Doob [11], who showed that posterior consistency occurs for all $\theta$ in a set of prior measure one. In this paper we consider a study of the non i.i.d. case based on Doob's result, specifically, the simple linear regression model. The martingale techniques are not applied here and we discuss the connection of posterior consistency with orthogonality of product measures.

Consistency is just the beginning of Bayesian asymptotics. Issues such as rates of convergence and asymptotic normality have received quite a lot of attention. Yet it appears that even at the level of consistency there are still issues that need to be clarified. In this paper, we review some conceptual issues in consistency of posterior distributions, and discuss different approaches to posterior consistency that have been investigated in the literature; we view this as a followup of Ghosal et al. [14]. We have attempted to elucidate those sufficient conditions to establish posterior consistency and tie up some loose ends on diverse conceptual issues in consistency of posterior distributions. The paper also contains some new results along with a brief commentary to the subject. In general, detailed proofs are omitted and given only when they are different from standard published materials or when the result is unpublished.

Section 2 contains a summary of some background material, some of the notations and assumptions used in the paper. Section 3 largely describes known results although some of the proofs are reorganized. The criteria based on the uniform strong law of large numbers is new and so far we are not aware of any significant application. The result might still be of interest because of its similarities to results on the consistency of nonparametric maximum likelihood estimates (NPMLEs) and also because it affords a natural extension to non i.i.d. cases. Section 4 specializes the results in Section 2 to the context of consistency. After a brief discussion of Schwartz's result, we discuss the known conditions for $L_1$ consistency and the relationship between these. Section 5 extends the Schwartz theorem to improper priors and formal posteriors. The result is new even in the parametric case. We have not pursued conditions for stronger consistency because improper priors usually arise in finite dimensional situations where weak and strong consistency coincide. Section 6 contains an example. All the general consistency results in the literature actually establish exponential consistency. In Section 6 we give an example where consistency obtains but not exponential consistency. The example surprised us as we had believed that, at least in the i.i.d. case, consistency would always be at an exponential rate.

In the last section we study the extension of consistency results to a non i.i.d. case. We give an example to show that the analogue of Doob's theorem will not always hold, and we prove a Doob theorem for the linear regression model with nonparametric errors. We also briefly discuss an extension of the theorem of Walker.

## 2. Preliminaries

In the setup that we consider, $\Theta$ is the parameter space; $\{f_\theta : \theta \in \Theta\}$ is a family of densities with respect to a $\sigma$-finite measure $\mu$ on a measurable space $\mathcal{X}$. We will use $P_\theta$ to denote the probability distribution generated by $f_\theta$. Throughout the paper we assume that $\Theta$ and $\mathcal{X}$ are complete separable metric spaces and we also assume that $\theta \mapsto f_\theta$ is 1-1 and $(\theta, x) \mapsto f_\theta(x)$ is measurable.

The affinity, $\text{Aff}(f,g)$, between any two densities is defined as $\text{Aff}(f,g) = \int \sqrt{fg} d\mu$. Let $\Pi$ be a prior distribution, i.e., a probability measure on $\Theta$. Given $\theta, X_1, X_2, \ldots, X_n$ are assumed to be i.i.d. $P_\theta$. $f_\theta^{(n)}(x_1, x_2, \ldots, x_n)$ will stand for the joint density $\prod_{i=1}^n f_\theta(x_i)$.

The Kullback–Leibler (KL) divergence is denoted by $K(\theta_0, \theta) = E_{\theta_0} \log(f_\theta/f_{\theta_0})$. A KL neighborhood $K_\epsilon(\theta_0)$ of $\theta_0$ is denoted by $\{\theta : K(\theta_0, \theta) < \epsilon\}$.

**Definition 2.1.** A point $\theta_0$ is said to be in the KL support of $\Pi$ if for all $\epsilon > 0, \Pi(K_\epsilon(\theta_0)) > 0$.

The posterior distribution $\Pi(A|X_1, X_2, \ldots, X_n)$, the version that we consider, is given by the following. For any measurable subset $A$ of $\Theta$,

$$(2.1) \qquad \Pi(A|X_1, X_2, \ldots, X_n) = \frac{J_A(X_1, X_2, \ldots, X_n)}{J(X_1, X_2, \ldots, X_n)}$$

where

$$J_A(X_1, X_2, \ldots, X_n) = \int_A \frac{f_\theta^{(n)}}{f_{\theta_0}^{(n)}}(X_1, X_2, \ldots, X_n) \Pi(d\theta)$$

and
$$J(X_1, X_2, \ldots, X_n) = \int_\Theta \frac{f_\theta^{(n)}}{f_{\theta_0}^{(n)}}(X_1, X_2, \ldots, X_n) \Pi(d\theta).$$

## 3. Exponential decrease to 0

We begin with a review of results that provide conditions under which, for a measurable subset $A$ of $\Theta$, $\Pi(A|X_1, X_2, \ldots, X_n)$ goes to 0 exponentially with $P_{\theta_0}^\infty$ probability 1.

**Definition 3.1.** Let $\theta_0 \in \Theta$ and let $P_{\theta_0}^\infty$ stand for the joint distribution of $\{X_i\}_{i=1}^\infty$ when $\theta_0$ is the true value of $\theta$. Then $\Pi(A|X_1, X_2, \ldots, X_n)$ is said to go to 0 exponentially with $P_{\theta_0}^\infty$ probability 1, if there exists a $\beta > 0$ such that
$$P_{\theta_0}^\infty \left( \{\Pi(A|X_1, X_2, \ldots, X_n) > e^{-n\beta} \text{ i.o. }\} \right) = 0$$
where *i.o.* stands for "infinitely often."

Proposition 3.2 goes back to [12] and [25]. For a proof see [19, Lemma 4.4.1].

**Proposition 3.2.** *If $\theta_0$ is in the KL support of $\Pi$ then for all $\beta > 0$,*
$$\lim_{n \to \infty} e^{n\beta} J(X_1, X_2, \ldots, X_n) = \infty \text{ a.s. } P_{\theta_0}^\infty.$$

Proposition 3.2 shows that the Kullback–Leibler support condition takes care of the denominator in (2.1). The exponential convergence to 0 would follow if it can be established that there exists $\beta_0 > 0$ such that $e^{n\beta_0} J_A(X_1, X_2, \ldots, X_n) \to 0$ a.s. $P_{\theta_0}^\infty$. We explore sufficient conditions to achieve this.

**Definition 3.3.** For a probability measure $\nu$ on $\theta$, let $q_\nu^{(n)}$ be the marginal density of $X_1, \ldots, X_n$,
$$q_\nu^{(n)}(x_1, x_2, \ldots, x_n) = \int_\Theta f_\theta^{(n)}(x_1, x_2, \ldots, x_n) \nu(d\theta).$$

**Definition 3.4.** Let $A \subset \Theta$ and $\delta > 0$. The set $A$ and $\theta_0$ are said to be *strongly $\delta$ separated* if for any probability $\nu$ on $A$,
$$\text{Aff}(f_{\theta_0}, q_\nu^{(1)}) < \delta.$$

The relationship $H^2(f, g) = 1 - 2\text{Aff}(f, g)$ between the Hellinger distance $H(f, g)$ and the Affinity $\text{Aff}(f, g)$ shows that $\text{Aff}(f_{\theta_0}, q_\nu^{(1)}) < \delta$ is equivalent to $H^2(f_{\theta_0}, q_\nu^{(1)}) > 1 - \delta$. Say that $A$ and $\theta_0$ are strongly separated if they are strongly $\delta$ separated for some $\delta > 0$.

**Example 3.5.** Suppose that the $L_1$ distance between $f_{\theta^*}$ and $f_{\theta_0}$ is larger than $\delta^*$ for some $\delta^* > 0$, $\|f_{\theta^*} - f_{\theta_0}\| > \delta^*$. Let
$$A = \left\{ \theta : \|f_{\theta^*} - f_\theta\| < \frac{\delta^*}{2} \right\}.$$

It is easy to see that $A$ is strongly separated from $\theta_0$ for every $\nu$ on $A$.

We begin by isolating a useful consequence of strong separation. The underlying idea is in [32]. Note that the argument is essentially analytic and does not use Hoeffding's inequality as in [19]. Lemma 3.6, we believe, can be extended to non-i.i.d and even to non-independent cases. We do not pursue this here but will briefly return to it in Section 7.

**Lemma 3.6.** *If $\theta_0$ and $A$ are strongly $\delta$ separated then for all probability $\nu$ on $A$, for all $n$,*

$$\text{Aff}(f_{\theta_0}^{(n)}, q_\nu^{(n)}) < e^{-n\beta_0}, \text{ where } \beta_0 = -\log \delta. \tag{3.1}$$

*Proof.* The proof is straightforward by induction on $n$, combined with the definition of strong separation. □

**Remark 3.1.** The conclusion of Lemma 3.6 holds with $\beta_0 = -\log \delta/k$ if for all $\nu$, for some $k$, $\text{Aff}(f_{\theta_0}^{(k)}, q_\nu^{(k)}) < \delta$, i.e., $A$ and $\theta_0$ are strongly separated for the parametrization $\theta \mapsto f^{(k)}$.

The next result is the celebrated result of Schwartz [25] stated in terms of strong separation. A result of LeCam [21] shows that it is equivalent to the formulation of Schwartz involving an unbiased test for testing $H_0 : \theta = \theta_0$ vs. $H_1 : \theta \in A$. LeCam's theorem is proved using the Hahn–Banach theorem so is essentially an existence result. Hence the point of view of strong separation could be an easier condition to verify in some situations.

**Theorem 3.7** (Schwartz). *If*

(1) *$\theta_0$ is in the KL support of $\Pi$,*
(2) *for some $k$, $A$ and $\theta_0$ are strongly separated for the parametrization $\theta \mapsto f^{(k)}$.*

*Then $\Pi(A|X_1, X_2, \ldots, X_n)$ goes to 0 exponentially a.e. $P_{\theta_0}^\infty$.*

*Proof.* Let $\Pi^*$ be the probability measure obtained by restricting $\Pi$ to $A$ and normalizing it. Then

$$\begin{aligned} P_{\theta_0}(\sqrt{J_A} > e^{-n\gamma}) &\leq e^{n\gamma} E_{\theta_0}(\sqrt{J_A}) \\ &= e^{n\gamma} \sqrt{\Pi(A)} \text{Aff}(f_{\theta_0}^{(n)}, q_{\Pi^*}^{(n)}) \\ &\leq \sqrt{\Pi(A)} e^{n\gamma} e^{-n\beta_0}. \end{aligned}$$

Taking $\gamma = \beta_0/4$, it follows easily that

$$P_{\theta_0}(\sqrt{J_A} > e^{-n\gamma} \text{ i.o.}) = 0.$$

The proof can be completed easily using Proposition 3.2. For details see [19]. □

Proposition 3.2 and Lemma 3.6 easily give the following theorem of Walker [32].

**Theorem 3.8.** *If*

(1) *$\theta_0$ is in the KL support of $\Pi$.*
(2) *If $A = \cup_{i \geq 1} A_i$ such that:*

   (a) *For some $\delta > 0$ all the $A_i$'s are strongly $\delta$ separated from $\theta_0$ and*
   (b) *$\sum_{i \geq 1} \sqrt{\Pi(A_i)} < \infty$.*

*Then $\Pi(A|X_1, X_2, \ldots, X_n)$ goes to 0 exponentially a.e. $P_{\theta_0}^\infty$.*

*Proof.* It follows by noting

$$P_{\theta_0}(\sqrt{J_A} > e^{-n\gamma}) \leq e^{n\gamma} E_{\theta_0}(\sqrt{J_A})$$

$$\leq e^{n\gamma} E_{\theta_0}\left(\sqrt{\sum_i J_{A_i}}\right) \leq e^{n\gamma} \sum_i E_{\theta_0}\left(\sqrt{J_{A_i}}\right)$$

(3.2)
$$= e^{n\gamma} \sum_i \sqrt{\Pi(A_i)} \mathrm{Aff}(f_{\theta_0}^{(n)}, q_{\Pi_i^*}^{(n)})$$

$$\leq e^{n\gamma} e^{-n\beta_0} \sum_i \sqrt{\Pi(A_i)},$$

where $\Pi^*$ in (3.2) is the normalized restriction of $\Pi^*$ to $A_i$. □

The next theorem gives another set of sufficient conditions, in terms of the uniform Strong Law of Large Numbers (SLLN), for the posterior probability of a set to go to 0 exponentially. The conditions are stronger than those of Schwartz [25]. They are similar in spirit to the conditions used in the study of Hellinger consistency of NPMLEs (see [30]) and suggest a parallel between consistency of NPMLEs and posterior consistency.

**Theorem 3.9.** *Let $A \subset \Theta$. If*

(1) *$\theta_0$ is in the KL support of $\Pi$.*
(2) *$\mathrm{Aff}(f_{\theta_0}, f_\theta) < \delta$ for all $\theta \in A$.*
(3) $\sup_{\theta \in A} \left| \int \sqrt{\frac{f_\theta}{f_{\theta_0}}}(x) dP_n - \mathrm{Aff}(f_{\theta_0}, f_\theta) \right| \to 0$ *a.s $P_{\theta_0}^\infty$, where $P_n$ is the empirical distribution obtained from $X_1, X_2, \ldots, X_n$.*

*Then $\Pi(A|X_1, X_2, \ldots, X_n)$ goes to 0 exponentially a.e. $P_{\theta_0}^\infty$.*

*Proof.* Note that $\int g(x) dP_n = (1/n) \sum_{i=1}^n g(X_i)$ for arbitrary function $g(x)$. Thus,

$$J_A = \int_A \prod_1^n \frac{f_\theta}{f_{\theta_0}}(X_i) \Pi(d\theta)$$

$$= \int_A \exp\left\{2n \int \log \sqrt{\frac{f_\theta}{f_{\theta_0}}}(x) dP_n \right\} \Pi(d\theta)$$

since $\log x \leq x - 1$

$$\leq \int_A \exp\left\{2n \int \left(\sqrt{\frac{f_\theta}{f_{\theta_0}}}(x) - 1\right) dP_n \right\} \Pi(d\theta).$$

Take $\delta^* = 1 - \delta$. By assumptions 2 and 3, for all large $n$,

$$\sup_{\theta \in A} \sqrt{\frac{f_\theta}{f_{\theta_0}}} dP_n \leq \sup_{\theta \in A} \left\{ \left| \sqrt{\frac{f_\theta}{f_{\theta_0}}} dP_n - \mathrm{Aff}(f_{\theta_0}, f_\theta) \right| + \mathrm{Aff}(f_{\theta_0}, f_\theta) \right\}$$

$$\leq \frac{\delta^*}{2} + 1 - \delta^* = 1 - \delta^*/2,$$

which in turn implies that $J_A < \Pi(A) \exp(-n\delta^*/2)$. □

**Proposition 3.10.** *Conditions (2) and (3) of Theorem 3.9 imply that there exists a uniformly consistent test for $H_0 : \theta = \theta_0$ vs. $H_1 : \theta \in A$.*

*Proof.* Choose $\delta_0$ such that $\delta + \delta_0 = \beta_0 < 1$.

Let

$$C = \left\{(x_1, x_2, \ldots, x_n) : \sup_{\theta \in A} \left| \int \sqrt{\frac{f_\theta}{f_{\theta_0}}}(x) dP_n - \text{Aff}(f_{\theta_0}, f_\theta) \right| < \delta_0 \right\}.$$

By assumption (3) of Theorem 3.9, for any $\epsilon > 0$ and sufficiently large $n$, $P_{\theta_0}(C) > 1 - \epsilon$. For each $(x_1, x_2, \ldots, x_n)$ in $C$, for all $\theta \in A$,

$$\frac{1}{n} \sum \sqrt{\frac{f_\theta}{f_{\theta_0}}}(x_i) \leq \sup_{\theta \in A} \text{Aff}(f_{\theta_0}, f_\theta) + \delta_0$$

so that

$$\sum \left( \sqrt{\frac{f_\theta}{f_{\theta_0}}}(x_i) - 1 \right) < n(\delta + \delta_0 - 1) = -n\beta_0.$$

Therefore, for $\theta \in A$,

$$P_\theta(C) = \int_C \frac{f_\theta^{(n)}}{f_{\theta_0}^{(n)}}(x_i) f_{\theta_0}^{(n)}(x_i) \prod \mu(dx_i) < P_{\theta_0}(C) e^{-2n\beta_0}. \qquad \square$$

**Remark 3.2.** Salinetti [24] has used the notion of hypo convergence to study consistency of posterior and consistency of maximum likelihood estimates. A somewhat related result is due to Ghosal and van der Vaart who show that her condition is related to Schwartz's testing condition in the discussion of [24].

While the results discussed so far deal with sufficient conditions for $\Pi(A|X_1, \ldots, X_n)$ to go to 0 exponentially, the next basic result due to Barron [2] gives conditions that are both necessary and sufficient.

**Theorem 3.11** (Barron). *$A \subset \Theta$. Assume that $\theta_0$ is in the KL support of $\Pi$. Then the following are equivalent.*

(i) *There exists a $\beta_0$ such that*

$$P_{\theta_0}\{\Pi(A|X_1, X_2, \ldots, X_n) > e^{-n\beta_0} \ i.o.\} = 0.$$

(ii) *There exist subsets $V_n, W_n$ of $\Theta$, positive numbers $c_1, c_2, \beta_1, \beta_2$, and a sequence of tests $\{\phi_n\}$, $\phi_n$ based on $n$ observations, such that*

(a) $A \subset V_n \cup W_n$,
(b) $\Pi(W_n) \leq C_1 e^{-n\beta_1}$, and
(c) $P_{\theta_0}\{\phi_n > 0 \ i.o.\} = 0$ and $\inf_{f \in V_n} E_f \phi_n \geq 1 - c_2 e^{-n\beta_2}$.

## 4. Consistency

As before, $\Pi$ stands for a prior and $\{\Pi(\cdot|X_1, X_2, \ldots, X_n)\}$ denotes a sequence of posterior distributions. The sequence of posteriors is said to be consistent at $\theta_0$ if $\{\Pi(U|X_1, X_2, \ldots, X_n)\} \to 1$ a.s.$P_{\theta_0}^\infty$ for all neighborhoods $U$ of $\theta_0$.

Typically the parametrization $\theta \mapsto f_\theta$ turns out to be continuous when the space of densities is endowed with weak convergence or with the $L_1$ or the Hellinger metric. Consequently the neighborhoods of interest are those that arise from weak or $L_1$ topology.

In view of the last section, what is required then is to verify that the conditions developed in the last section apply to neighborhoods.

Let $g(x)$ be a bounded measurable function and define

$$(4.1) \qquad A_g = \left\{ \theta : \int g(x) f_\theta(x) \mu(dx) - \int g(x) f_{\theta_0}(x) \mu(dx) \geq \epsilon \right\}.$$

Clearly $A$ is strongly $\epsilon$ separated from $\theta_0$ and hence if $\theta_0$ is in the KL support of $\Pi$ then by Theorem 3.7 the posterior probability of $A$ goes to 0 exponentially. If $U$ is a weak neighborhood then $U^c$ is a finite union of sets of the type displayed in (4.1). This establishes exponential consistency for weak neighborhoods.

Consider the $L_1$ neighborhood

$$U = \{\theta : \|f_\theta - f_{\theta_0}\| < \epsilon\}.$$

In this case, in general, $U^c$ cannot be expressed as a finite union of sets strongly separated from $\theta_0$. Unlike the case of weak neighborhoods, in this case we need conditions beyond requiring that $\theta_0$ is in the KL support of $\Pi$.

Theorem 3.8 can be easily adapted in this context.

**Theorem 4.1.** *Assume*

(1) *$\theta_0$ is in the KL support of $\Pi$.*
(2) *For all $\delta > 0$, there exist sets $A_1, A_2, \ldots$ such that the diameter of $A_i$, $\text{diam}(A_i) < \delta$, $\bigcup A_i = \Theta$ and $\sum \sqrt{\Pi(A_i)} < \infty$.*

*Then for any $L_1$ neighborhood $U$ of $\theta_0$, the posterior probability of $U^c$ goes to 0 exponentially a.e. $P_{\theta_0}^\infty$.*

The theorem follows from observing that if $U$ is an $\epsilon$ neighborhood, then taking $\{A_i\}_{i=1}^\infty$ for $\delta = \epsilon/3$ it is easily seen that the $A_i$'s that have non-empty intersection with $U^c$ cover $U^c$. These $A_i$'s satisfy the assumptions of Theorem 3.8.

**Definition 4.2** (Bracketing entropy). Let $\Gamma \subset \Theta$. For a $\delta > 0$ define the bracketing entropy $\mathcal{H}(\Gamma, \delta)$ to be the logarithm of the minimum integer $k$ such that, there exist non negative functions $f_1^U, f_2^U, \ldots, f_k^U$ satisfying

(1) $\int f_i^U(x) \mu(dx) < 1 + \delta$,
(2) for each $\theta$ there exists $i$ such that $f_\theta \leq f_i^U$.

**Definition 4.3** (Metric entropy). Let $\Gamma \subset \Theta$. For $\delta > 0$ the Metric entropy $J(\Gamma, \delta)$ is defined to be the logarithm of minimum of all integers $k$ such that there exist densities $f_1^*, f_2^*, \cdots, f_k^*$ such that for each $\theta$ there exists $i$ such that $\|f_\theta - f_i^*\| < \delta$.

If $\theta_0$ is in the KL support of $\Pi$ then each of the three conditions listed below ensures that the posterior is exponentially $L_1$ consistent. The first condition (W) is from Walker's theorem, Theorem 3.8, the next (BSW) is due to [3] and the third (GGR) appears in [15]. A formal statement and proof can be found in [3] and [15].

(W) For each $\delta > 0$, there exist sets $A_1, A_2, \ldots$ such that $\bigcup A_i = \Theta$, $L_1$-diameter of $\{f_\theta : \theta \in A_i\} < \delta$ and $\sum_i \sqrt{\Pi(A_i)} < \infty$.

(BSW) For each $\epsilon > 0$, there exist $\Theta_n \subset \Theta$, and $C, c_1, c_2, \delta$ all positive such that

(a) $\Pi(\Theta_n^c) < e^{-nc_2}$.
(b) $\mathcal{H}(\Theta_n, \delta) \leq nc$ for $c < ([\epsilon - \sqrt{\delta}]^2 - \delta)/2, \delta < \epsilon^2/4$.

(GGR) If for each $\epsilon > 0$, there is a $0 < \delta < \epsilon, c_1, c_2, \beta < \epsilon^2/2$ and $\Theta_n$ such that

(a) $\Pi(\Theta_n^c) < c_1 e^{-n\beta}$.
(b) $J(\Theta_n, \delta) \leq n\beta$.

The next theorem shows that both (W) and (BSW) imply (GGR). A proof that (W) $\Rightarrow$ (GGR) was also communicated to us by Ghosal, S. [personal communication].

**Theorem 4.4.** *(W)$\Rightarrow$(GGR) and (BSW) $\Rightarrow$ (GGR).*

*Proof.* (W)$\Rightarrow$(GGR).

Assume without loss of generality that $\Pi(A_i) = \Pi_i$ is decreasing in $i$ and let $\sum \sqrt{\Pi_i} = c < \infty$. Set

$$\Theta_n = \bigcup_1^{k_n} A_i.$$

Since the $L_1$-diameter of $\{f_\theta : \theta \in A_i\} < \delta$, it is easy to see that $J(\Theta_n, \delta) < \log k_n$. Thus, by taking $k_n = \exp(n\beta)$ one then obtains sieves with the properties required by (GGR).

Next, we shall argue that $\Pi(\Theta_n^c) = \Pi\left(\bigcup_{i > k_n} A_i\right) \leq 2c^2/k_n$. Note that, for any $j$, $j\sqrt{\Pi_j} \leq \sum_{i=1}^j \sqrt{\Pi_i} \leq c$ so that $j \leq c/\sqrt{\Pi_j}$. Therefore,

$$\Pi\left(\bigcup_{j > k_n} A_i\right) \leq \sum_{j > k_n} \Pi_j \leq c^2 \sum_{j > k_n} \frac{1}{j^2} \leq \frac{2c^2}{k_n}.$$

(BSW)$\Rightarrow$(GGR)

Let $f_1, f_2, \ldots, f_k$ be functions such that $\int f_i = 1 + c_i < (1 + \delta)$ and such that for any $\theta \in \Gamma, \exists\ i$ such that $f_\theta \leq f_i$. Let $f_i^* = f_i/(1 + c_i)$. Then

$$\begin{aligned}
\|f_i^* - f_\theta\| &\leq \frac{1}{1 + c_i} \|f_i - (1 + c_i)f_\theta\| \\
&\leq \|f_i - f_\theta\| + c_i \\
&\leq 2\delta.
\end{aligned}$$

Hence $f_1^*, f_2^*, \ldots, f_k^*$ forms a $2\delta$ net for $\Gamma$ and $J(\Gamma, 2\delta) \leq \mathcal{H}(\Gamma, \delta)$. □

## 5. Improper priors and formal posteriors

Suppose that $\Pi$ is an improper prior on $\Theta$, that is, a $\sigma$-finite measure with $\Pi(\Theta) = \infty$. A formal posterior density given $X_1 = x_1, X_2 = x_2, \ldots X_n = x_n$ is defined as in Equation (2.1). This is of course well defined only if

$$J(x_1, x_2, \ldots, x_n) = \int_\Theta \frac{f_\theta^{(n)}}{f_{\theta_0}^{(n)}}(x_1, x_2, \ldots, x_n)\Pi(d\theta) < \infty.$$

This situation occurs widely in the context of noninformative priors (see for example, Ghosh and Ramamoorthi [18] and Kass and Wasserman [20]).

The next theorem shows that if $P_0$ is in the KL support of $\Pi$ then the posterior is weakly consistent. Improper priors largely arise in the context of finite dimensional regular models and in these situations weak consistency and strong consistency coincide. Hence, we do not develop conditions akin to (W), (BSW) or (GGR) for improper priors. First, Lemma 5.1 states a result of the KL support of $\Pi(\cdot|x)$.

**Lemma 5.1.** *Let $P_0$ is in the KL support of $\Pi$. Denote by $A = \{x : J(x) = \int f_\theta(x)\Pi(d\theta) < \infty\}$. Then, for $P_0$ almost all $x$ in $A$, $\theta_0$ is in the KL support of $\Pi(\cdot|x)$.*

*Proof.* Fix $\epsilon > 0$. Consider $E = \{x \in A : \Pi(K_\epsilon|x) = 0\}$. We shall show that $P_{\theta_0}(E) = 0$. Note that

$$\Pi(K_\epsilon|x) = \frac{\int_\Theta I_{K_\epsilon}(\theta)f_\theta(x)\Pi(d\theta)}{\int f_\theta(x)\Pi(d\theta)}.$$

Denoting by $\Pi^*$ the measure $\Pi(\cdot \cap K_\epsilon)/\Pi(K_\epsilon)$, since, for $x \in E$, $\Pi(K_\epsilon|x) = 0$, we have that

$$\Pi^*\{\theta : f_\theta(x) = 0\} = 1.$$

Consequently $\int_E \int_{K_\epsilon} f_\theta(y)\Pi^*(d\theta)(E)d\mu(y) = 0$. Interchanging the integrals, $\int_{K_\epsilon}[\int_E f_\theta(y)d\mu(y)]\Pi^*(d\theta) = 0$ and hence there exists some $\theta'$ such that $\int_E f_{\theta'}(y)d\mu(y) = 0$ so that $P_{\theta'}(E) = 0$. For every $\theta$ in $K_\epsilon$ $P_\theta$ dominates $P_{\theta_0}$, so $P_{\theta_0}(E) = 0$. Letting $\epsilon$ run through rationals, the lemma is established. □

**Theorem 5.2.** *Let $\Pi$ be an improper prior on $\Theta$. $\{f_\theta : \theta \in \Theta\}$ is a family of densities. Assume that the formal posterior is defined with $P_0^\infty$ probability one. Formally, if*

$$A_n = \{x_1, x_2, \ldots, x_n : J(x_1, x_2, \ldots, x_n) < \infty\} \text{ then } P_0^\infty(\bigcup A_n) = 1.$$

*If $\theta_0$ is in the KL support of $\Pi$ then the formal posterior is weakly consistent at $\theta_0$.*

*Proof.* By Lemma 5.1, for each $n$, except for those in a set of $P_{\theta_0}$ measure 0, for all $(x_1, x_2, \ldots, x_n) \in A_n$, $\theta_0$ is in the KL support of $\Pi(\cdot|(x_1, x_2, \ldots, x_n))$.

Since on $A_n$, $\Pi(\cdot|(x_1, x_2, \ldots, x_{n+1})) = \Pi_{(x_1,x_2,\ldots,x_n)}(\cdot|x_{n+1})$, the result follows. □

## 6. Example

All the results discussed so far are related to exponential consistency. The next example shows that, even in the context of i.i.d. observations, the posterior can be consistent at a non-exponential rate.

Consider an example where we have a prior $\Pi$, $f_0$ is in the KL support of $\Pi$ and the posterior is not $L_1$ consistent, i.e., there is a set $A$ which is a complement of a neighborhood of $f_0$ and whose posterior does not go to 0. Such an example appears, for instance, in Barron et al. [3].

Consider the prior to $\Pi^* = .5\delta_{f_0} + .5\Pi$. Then by Doob's theorem the posterior of $A$ goes to 0. It cannot go exponentially, for if it does, by Barron's theorem (e.g. [2] and [19, Theorem 4.4.3]), there would be sieves $V_n$ and sets $U_n$ of exponentially small $\Pi^*$ probability that cover $A$. These properties also carry over to $\Pi$ and now the first part of Barron's result would imply that the original prior $\Pi$ is itself consistent, in fact, exponentially consistent.

## 7. Independent but non-identically distributed models

### 7.1. Extension to posterior consistency

Here we look at the setup where, as before, $\Theta$ is a parameter space and $\Pi$ is a prior on $\Theta$. Given $\theta$, we assume that $X_1, X_2, \ldots$ are independent with $X_i$ distributed as $f_{i,\theta}$.

All the results discussed so far can be easily adapted, but not necessarily easily applied in the non-identically distributed case. As in Section 1, the posterior can be written as the ratio of two integrals. A stronger form of KL support (for instance, see Choudhuri et al. [9] and Amewou-Atisso et al. [1]) takes care of the denominator. It is not clear if there is a simple version of the (GGR) type of sufficient condition. Instead, those results for independent but non-identically distributed models as in Amewou-Atisso et al. [9], Choudhuri et al. [1], Ghosal and Roy [16] and Choi and Schervish [8], tried to establish the existence of uniformly consistent tests directly, which makes the numerator in the ratio of two integrals decrease to 0 exponentially.

Alternatively, LeCam [22] and Birge [4] showed that for independent non-identically distributed variables, tests with exponentially small errors exist when we use the average squared Hellinger distance to separate densities and convex sets. That is, uniformly consistent tests are always obtained if the entropy with such a distance is controlled. In the recent paper by Ghosal and van der Vaart [17], (GGR) type results have been investigated in the test construction for the convergence rates of posterior distributions for non i.i.d. observations.

On the other hand, Walker's sufficient conditions are easily adaptable in this case. Note that the proof of Lemma 3.6 does not require the assumption of the identically distributed observations; hence Theorem 3.8 easily follows to this case. We state it formally below.

**Theorem 7.1.** *If $A = \bigcup_{i \geq 1} A_i$ such that*

1. *For some $\delta > 0$ all the $A_i$'s are strongly $\delta$ separated from $\theta_0$ for the model $\theta \mapsto f_{i,\theta}$ and*
2. $\sum_{i \geq 1} \sqrt{\Pi(A_i)} < \infty$.

*Then for some $\beta_0 > 0$,*

$$e^{n\beta_0} \int_A \prod_{i=1}^n \frac{f_{i,\theta}(x_i)}{f_{i,\theta_0}(x_i)} \Pi(d\theta) \to 0 \ a.s \ \prod_{i=1}^\infty P_{i,\theta_0}.$$

Similar results to Theorem 7.1 along with regression problems have been discussed in Walker [33].

**Example 7.2** (Orthogonal series expansion). Let

(7.1) $$Y_i = \eta(X_i) + \epsilon_i, \ i = 1, \ldots, n$$

where the $\epsilon_i$'s are assumed to be independent $N(0,1)$ random variables, the $X_i$'s are sampled from a known probability distribution, and $\eta(\cdot)$ is a regression function. An orthogonal series expansion for the regression function $\eta(x)$ is a representation of $\eta(x)$ by an infinite sum,

$$\eta(x) = \sum_{j=1}^\infty \eta_j \phi_j(x),$$

where $\{\phi_j(x)\}_{j=1}^{\infty}$ is an orthonormal basis for an $L^2$ space containing $\eta$. Regarding either posterior consistency or rate of convergence of posterior distributions, this model has been investigated by Shen and Wasserman [26], Walker [33] and Choi and Schervish [8].

Let $\{\phi_j(\cdot)\}_{j=1}^{\infty}$ be an orthonormal basis for $L^2[0,1]$ such that for some $C > 0$, $\sup_{x \in [0,1]} |\phi_j(x)| \leq C$ for all $j$.

In this case, we consider the following $\delta$-covering of $\Omega$, a union of sets of the type

(7.2) $$\{\psi \ : \ n_j \delta_j < \psi_j < (n_j + 1)\delta_j, \ j = 1, 2, \ldots\},$$

which was also examined for Hellinger consistency in density estimation problems from infinite-dimensional exponential families in Walker [32] and regression problems in Walker [33]. Based on (7.2), the condition (b) in Theorem 7.1 can be verified as in Section 6.1 [32]. When the regression function is uniformly bounded, the $L_1$(or Hellinger) neighborhood of the true density $f_{\theta_0}$ becomes equivalent to the $L_1$ neighborhood of the true regression function $\eta_0$. Therefore, by considering a $\delta$-covering in (7.2) and its corresponding prior probability, two conditions (a) and (b) are easily verified. Hence, the conclusion of Theorem 7.1 is achieved when $A$ is in the $L_1$ neighborhood of the true density generating the regression model (7.1).

**Example 7.3** (Gaussian process regression). Gaussian process regression is one of the popular approaches to Bayesian nonparametric regression problems, and it is used to model the regression function $\eta(x)$ as a Gaussian process a priori. Posterior consistency based on Gaussian processes has been established in Ghosal and Roy [16] and Choi [7] for nonparametric binary regression, Tokdar and Ghosh [29] for density estimation and Choi [6] and Choi and Schervish [8] for nonparametric regression. Interestingly, all the results mentioned above have been based on constructing uniformly consistent tests rather than the condition (b) in Theorem 7.1. The challenges in the study of posterior consistency based on Gaussian processes is to find a rate that a prior probability shrinks as we consider a sequence of $\delta$-coverings that satisfies the condition (b). In this case, the important task to be achieved is obtaining the exponentially small lower bound for small balls of Gaussian processes. There is a recent investigation in this regard (e.g. see Li and Shao [23] and van der Vaart and van Zanten [31]). It would be interesting to explore if this difficulty in verifying (b) under Gaussian process priors can be bypassed when we apply Theorem 7.1.

### 7.2. Doob's theorem

Doob [11] showed that when $\Theta$ is the parameter space and given $\theta$, $\mathbf{X}_1, \mathbf{X}_2, \ldots$, are i.i.d. $P_\theta$ then, for any sequence of posterior distributions $\Pi(\cdot|\mathbf{X}_1, \mathbf{X}_2, \ldots, \mathbf{X}_n)$, under mild set theoretic assumptions (for instance when $\mathcal{X}$ and $\Theta$ are Borel subsets of Polish spaces) for any prior $\Pi$, there is a $\Theta_\Pi \subset \Theta$ with $\Pi(\Theta_\Pi) = 1$ such that the posterior is consistent at all $\theta \in \Theta_\Pi$. In what follows we explore the analogue of Doob's theorem in independent non-identically distributed models.

To change the notation a bit, given $\theta$ in $\Theta$, let $\mathbf{Y}_1, \mathbf{Y}_2, \ldots,$ be $\mathcal{Y}$ valued random variables with joint distribution $P_{\theta,\infty}$. For any prior $\Pi$ on $\Theta$, denote by $\lambda_\Pi$ the joint distribution induced on $\Theta \times \mathcal{Y}^\infty$ by $\Pi$ and $\{P_{\theta,\infty} : \theta \in \Theta\}$. We will denote the elements of $\mathcal{Y}^\infty$ by $\mathbf{y}$ and of $(\mathbf{Y}_1, \mathbf{Y}_2, \ldots,)$ by $\mathbf{Y}$. As before $\Pi(\cdot|\mathbf{Y}_1, \mathbf{Y}_2, \ldots, \mathbf{Y}_n)$ will stand for a fixed version of the posterior distribution of $\theta$.

By going through an appropriate countable set of continuous functions $g$ and applying the martingale convergence theorem to each posterior mean of $g(\theta)$, it can be seen that there is a conditional probability $\Pi^*(\cdot|\boldsymbol{y})$ such that for all $\boldsymbol{y}$ outside a $\lambda_\Pi$ null set

$$\Pi(\cdot|\mathbf{Y}_1, \mathbf{Y}_2, \ldots, \mathbf{Y}_n) \stackrel{weakly}{\to} \Pi^*(\cdot|\mathbf{y}).$$

Clearly the posterior is consistent at $\theta$ if $\Pi^*(\cdot|\boldsymbol{y}) = \delta_\theta$ a.e. $P_{\theta,\infty}$.

**Proposition 7.4.** *Consider the following two sets of statements for a given prior $\Pi$:*

1. *There is a set $\Theta_\Pi$ with $\Pi(\Theta_\Pi) = 1$ and the posterior is consistent at all $\theta \in \Theta_\Pi$.*
2. *There is a set $\Theta_\Pi$ with $\Pi(\Theta_\Pi) = 1$ and a measurable set $E^\Pi \subset \Theta \times \mathcal{Y}^\infty$ such that*
   (a) *for each $\theta \in \Theta_\Pi, P_{\theta,\infty}(E_\theta^\Pi) = 1$,*
   (b) *$E_\theta^\Pi \cap E_{\theta'}^\Pi = \emptyset$ for $\theta \neq \theta'$.*

*The two sets of statements are equivalent.*

*Proof.* Suppose (1) holds. Then it is easy to verify that the set

$$E^\Pi = \{(\theta, \boldsymbol{y}) : \Pi^*(\cdot|\boldsymbol{y}) = \delta_\theta, \ (\theta \in \Theta_\Pi)\}$$

is measurable and satisfies the conditions in (2).

On the other hand if (2) holds then define $\phi(\boldsymbol{y}) = \theta$ if $(\theta, \boldsymbol{y}) \in E^\Pi$. Then, using a result from set theory [28, Theorem 4.5.7], it can be shown that $\phi$ is measurable. It is easy to verify that $\tilde{\Pi}$ defined by

$$\tilde{\Pi}(\cdot|\boldsymbol{y}) = \delta_{\phi(\boldsymbol{y})}$$

is a version of the conditional distribution of $\theta$ given $\mathbf{Y}$ and hence $\Pi^*(\cdot|\boldsymbol{y}) = \tilde{\Pi}(|\boldsymbol{y})$ a.e. $\lambda_\Pi$. An application of Fubini's theorem yields the result. □

Our interest is in establishing (1) for all priors $\Pi$ and it is convenient to work with a stronger version of (2) by seeking a decomposition that does not depend upon $\Pi$. Formally,

**Proposition 7.5.** *Let $\Pi$ be a prior for $\Theta$.*
*Suppose there exists a measurable set $E \subset \Theta \times \mathcal{Y}^\infty$ such that*

1. *For each $\theta \in \Theta, P_{\theta,\infty}(E_\theta) = 1$ where $E_\theta$ is the $\theta$-section $\{\boldsymbol{y} : (\theta, \boldsymbol{y}) \in E\}$.*
2. *$E_\theta \cap E_{\theta'} = \emptyset$ for $\theta \neq \theta'$.*

*Then there is a set $\Theta_\Pi$ with $\Pi$ measure 1, such that the posterior is consistent at all $\theta \in \Theta_\Pi$.*

Thus, Doob's theorem is intimately related to uniform orthogonality of $\{P_{\theta,\infty} : \theta \in \Theta\}$. There is a wide literature on singularity and mutual absolute continuity of measures on infinite product spaces ([27] and [5]). This literature in general deals with pairwise orthogonality whereas Proposition 7.5 requires uniform orthogonality. The step from pairwise to uniform orthogonality can be formidable. Yet we feel that some of these results are likely to be useful in establishing Doob-type theorems in the non-i.i.d. set up.

Motivated by Proposition 7.5, we present an example where the Doob-type theorem fails to hold. On the positive side, Proposition 7.5 enables us to prove a theorem for linear regression models with nonparametric errors.

The case that we consider is

$$Y_i = \alpha + \beta x_i + \epsilon_i, \qquad i = 1, 2, \ldots$$

where

1. $x_1, x_2, \ldots$, are fixed non-random design points.
2. $\epsilon_1, \epsilon_2, \ldots$ are independent and identically distributed random variables with a probability density symmetric around 0.

**Example 7.6.** Suppose $\sum_i x_i^2 < \infty$ and $\epsilon_i \sim N(0, 1)$, and let $\alpha = 0$. In this case it follows from a result of Shepp [27] that $\prod_1^\infty N(\beta x_i, 1)$ are mutually absolutely continuous. Hence the decomposition required by Proposition 7.4 fails and Doob's theorem cannot hold.

The last example we consider is semiparametric regression, the linear regression model where the distribution of the noise is assumed to be unknown and thus needs to be estimated. This example has been investigated in terms of posterior consistency, following from the generalization of the Schwartz theorem in Amewou-Atisso et al. [1]. We revisit this example in Theorem 7.7 and show that the Doob-type theorem holds with an assumption on the fixed non-random design points, similar to that of [1].

Let Assumption A be defined as the following: There exists $\epsilon_0 > 0$ such that the covariate values $x_i$'s satisfy

$$\sum_i I_{(-\infty, -\epsilon_0)}(x_i) = \infty \text{ and } \sum_i I_{(\epsilon_0, \infty)}(x_i) = \infty.$$

**Theorem 7.7.** *Consider the model*

$$Y_i = \alpha + \beta x_i + \epsilon_i, \qquad i = 1, 2, \ldots$$

*where*

1. $x_1, x_2, \ldots$, *are fixed nonrandom design points.*
2. $\epsilon_1, \epsilon_2, \ldots$ *are i.i.d. variables with an unknown distribution of which density $f$ is symmetric, continuous at 0 and $f(0) > 0$.*

*If Assumption A holds, then given any prior $\Pi$ for $(\alpha, \beta, f)$, there is a set $\Theta_\Pi$ of $\Pi$ measure 1 such that the posterior is consistent at all $(\alpha, \beta, f) \in \Theta_\Pi$.*

*Proof.* Let $\mathcal{F}$ be all densities $f$ on the real line which are symmetric, continuous at 0 and $f(0) > 0$. Formally, we have as the parameter space $\Theta = R \times R \times \mathcal{F}$ and given $(\alpha, \beta, f)$, the $Y_i$'s are independent with $Y_i \sim f_{\alpha+\beta x_i}$, where $f_{\alpha+\beta x_i}(y) = f(y - (\alpha + \beta x_i))$.

We now construct a decomposition satisfying the conditions of Proposition 7.5.

Let $N_1 = \{n_1, n_2, \ldots\}$ be the subsequence of all $i$ with $x_i > \epsilon$ and $M_1 = \{m_1, m_2, \ldots\}$ be the subsequence of all $i$ with $x_i < -\epsilon$. Let $t$ be a real number and define $A_t = (t, \infty)$.

Note that the unknown parameter $\theta$ is the triple $\theta \equiv (\alpha, \beta, f)$. Following notations in Proposition 7.5, let $\Pi$ be a prior on $\Theta$ and let $E^\Pi$ be the set of all $(\alpha, \beta, f, \boldsymbol{y})$ such that for any real number $t$,

1. $\lim_{n \to \infty} \frac{1}{n} \sum_1^n I_{A_t}(y_i - (\alpha + \beta x_i)) = P_f(A_t).$
2. $\lim_{k \to \infty} \frac{1}{k} \sum_1^k I_{A_t}(y_{n_i} - (\alpha + \beta x_{n_i})) = P_f(A_t).$

3. $\lim_{l \to \infty} \frac{1}{l} \sum_1^l I_{A_t}(y_{m_i} - (\alpha + \beta x_{m_i})) = P_f(A_t)$.

Since $N_1, M_1$ are fixed subsequences and since it is enough to work with $t$-rational, $E^\Pi$ is easily seen to be measurable.

Further, for each $(\alpha, \beta, f)$, $[\prod_1^\infty P_{\alpha+\beta x_i}](E^\Pi_{\alpha,\beta,f}) = 1$, where $E^\Pi_{\alpha,\beta,f}$ is the $(\alpha, \beta, f)$-section $\{\boldsymbol{y} : (\alpha, \beta, f, \boldsymbol{y}) \in E^\Pi\}$ for each $(\alpha, \beta, f)$ as defined in Proposition 7.5. This follows by noting that under $[\prod_1^\infty P_{\alpha+\beta x_i}]$, $Y_1 - (\alpha + \beta x_1), Y_2 - (\alpha + \beta x_2), \ldots$ are i.i.d. with common density $f$. An application of the law of large numbers proves the claim.

We next argue that if $(\alpha_1, \beta_1, f_1) \neq (\alpha_2, \beta_2, f_2)$ then $E^\Pi_{\alpha_1, \beta_1, f_1} \cap E^\Pi_{\alpha_2, \beta_2, f_2} = \emptyset$.

If $\alpha_1 = \alpha_2, \beta_1 = \beta_2$ and $f_1 \neq f_2$, and if $\boldsymbol{y} \in E^\Pi_{\alpha,\beta,f_1} \cap E^\Pi_{\alpha,\beta,f_2}$ then a contradiction is easily obtained by considering a $t$ for which $P_{f_1}(A_t) \neq P_{f_2}(A_t)$.

Now suppose that for some $\Delta > 0$, $\alpha_1 - \alpha_2 > \Delta$ and $\beta_1 - \beta_2 > \Delta$. Clearly for every $n_i \in N_1$, $(\beta_1 - \beta_2)x_{n_i} > \Delta \epsilon$. Choose $\eta$ such that $\eta < \Delta \epsilon$ and $\inf_{|x| < \eta} f_1(x) > C > 0$.

Since $f$ is symmetric and $\eta > 0$, $P_{f_1}(A_\eta) < 1/2$. We will get a contradiction by showing that if $\boldsymbol{y} \in E^\Pi_{\alpha_1,\beta_1,f_1} \cap E^\Pi_{\alpha_2,\beta_2,f_2}$, then $P_{f_1}(A_\eta) \geq 1/2$.

If $\boldsymbol{y} \in E^\Pi_{\alpha_1,\beta_1,f_1} \cap E^\Pi_{\alpha_2,\beta_2,f_2}$, then for all $t$,

$$\text{(7.3)} \qquad \frac{1}{k} \sum_1^k I_{A_t}(y_{n_i} - (\alpha_1 + \beta_1 x_{n_i})) \to P_{f_1}(A_t),$$

$$\text{(7.4)} \qquad \frac{1}{k} \sum_1^k I_{A_t}(y_{n_i} - (\alpha_2 + \beta_2 x_{n_i})) \to P_{f_2}(A_t),$$

$$\begin{aligned} \alpha_1 + \beta_1 x_{n_i} &= \alpha_2 + \beta_2 x_{n_i} + (\alpha_1 - \alpha_2) + (\beta_1 - \beta_2)x_{n_i} \\ &\geq \alpha_2 + \beta_2 x_{n_i} + \eta \end{aligned}$$

and hence

$$\begin{aligned} I_{A_t}(y_{n_i} - (\alpha_1 + \beta_1 x_{n_i})) &\geq I_{A_t}(y_{n_i} - (\alpha_2 + \beta_2 x_{n_i} + \eta)) \\ &= I_{A_{t-\eta}}(y_{n_i} - (\alpha_2 + \beta_2 x_{n_i})). \end{aligned}$$

In particular with $t = \eta$,

$$I_{A_\eta}(y_{n_i} - (\alpha_1 + \beta_1 x_{n_i})) \geq I_{(0,\infty)}(y_{n_i} - (\alpha_2 + \beta_2 x_{n_i})).$$

Consequently

$$\begin{aligned} \frac{1}{k} \sum_1^k I_{A_t}(y_{n_i} - (\alpha_1 + \beta_1 x_{n_i})) &\geq \frac{1}{k} \sum_1^k I_{(0,\infty)}(y_{n_i} - (\alpha_2 + \beta_2 x_{n_i})) \\ &\to P_{f_2}(0, \infty) = \frac{1}{2}. \end{aligned}$$

The case when $\alpha_1 - \alpha_2 < \Delta, \beta_1 - \beta_2 > \Delta$ can be handled by considering the subsequence $M_1$. Similarly, the other remaining cases follow. □

**Acknowledgments.** This paper, like all my research endeavours, has benefitted from J. K. Ghosh. Bayesian nonparametrics is but a small part of all that I have learned from him, and I am but one among many who owe their intellectual development to his influence. My association with J. K. Ghosh, JKG as I call him, began more than thirty five years ago in the form of student and teacher. This equation has remained constant but over the years, on top of it, has developed a friendship that I greatly cherish. This article is dedicated to JKG with admiration, appreciation, affection and... gratitude.

R. V. Ramamoorthi